\documentclass[12pt]{amsart}
\begin{document}

\title[MOTION BY WEIGHTED MEAN CURVATURE IS AFFINE INVARIANT]
{MOTION BY WEIGHTED MEAN CURVATURE\\IS AFFINE INVARIANT}

\author{Jean E. Taylor}
\address{Mathematics Department, Rutgers University, Piscataway,
NJ 08855}
\email{taylor@math.rutgers.edu}

\thanks{Research at MSRI is supported in part by NSF grant
DMS-9022140.}

\begin{abstract}
Suppose curves are moving by curvature in a plane, but one embeds the
plane in $R^3$ and looks at the plane  from an angle.  Then circles
shrinking to a round point would appear to be ellipses shrinking to an
``elliptical point,''  and the surface energy would appear to be anisotropic
as would the mobility.  The result of this paper is that if one uses the
apparent surface energy and the apparent mobility, then the motion by
weighted curvature with mobility in the apparent plane is the same 
as motion by curvature in the original
plane but then viewed from the angle.  This result
applies not only to the isotropic case but to 
arbitrary surface energy functions and mobilities in the plane, to surfaces in
3-space, and (in the case that the surface energy function is twice
differentiable) to the case of motion viewed through distorted lenses
(i.e., diffeomorphisms) as well.  This result is to be contrasted with
an earlier result [ST] which states that for area-preserving affine 
transformations of the plane where the energy and mobility are NOT also
transformed, motion by curvature to the power 1/3 (rather than 1) is invariant.
\end{abstract}

\maketitle

\section{Introduction} 

Suppose that $\Phi$ is a norm on $R^2$, except that it is
not necessarily even. That is, 
$\Phi(\lambda n) = \lambda \Phi(n)$ for every $\lambda \ge 0$,
$\Phi(n)= 0$ if and only if $n=0$, and $\Phi(n+p) \le
\Phi(n) + \Phi(p)$ for all $n$ and $p$ in $R^2$.
We regard $\Phi$ as a given convex ``surface'' (actually curve)
energy function for curves in the
plane; the energy of a rectifiable oriented curve $C$ is $\int_{x\in C}
\Phi(n_C(x)) \,dx$, where $n_C$ is the unit tangent $T_C$
rotated clockwise by  $90^\circ$. 
The weighted (mean) curvature $\kappa_\Phi$ of $C$ is
defined in section 2; it is the precise way of expressing the concept
of the rate of decrease of energy of $C$ with ``volume'' (in
$R^2$, area) swept out under deformations of $C$. 
Surface energy and weighted mean curvature are natural physical
concepts, arising from the fact that for solids whose atoms are in lattices,
interfaces with different normal directions relative to the lattice
are quite different in structure.  Various situations such as thin
films are often assumed to be effectively two-dimensional.

Suppose that $M$ is also a norm that is not necessarily even, and that
$C(t)$ is a family of curves in $R^2$
moving by weighted curvature $\kappa_\Phi$ with mobility $M$.
That is, for each $t$ and $x\in C(t)$, the normal velocity of  $C(t)$
at $x$ is $M(n_{C(t)}(x)) 
\kappa_\Phi(C(t),x)$.  Motion by weighted mean curvature with
mobility is likewise a standard physical concept, arising in
situations such as grain growth where the materials on each side
of the interface do not differ in their bulk energy.  See [TCH] for
examples and references.

Suppose $A$ is any affine map of the plane to itself; write
$A(x) = L(x) + b$, where $L$ is a general linear map and $b$ a
translation.  For each 
$q \in R^2$, define $\Phi_A(q) = \Phi(L^T q)$
and $M_A(q) = M(L^T q)$.   We show here that if
$C_A(t)$ is the family of curves moving by weighted mean curvature
$\kappa_{\Phi_A}$ and mobility $M_A$ and with $C_A(0) = A(C(0))$,
then $A(C(t)) = C_A(t)$ for every $t$. We also give the generalization
to higher dimensions.

This result can be interpreted in the following way.
Suppose you put the plane in which the curve is moving into $R^3$ and
look at it from an angle.  Then you would in fact observe
motion by weighted curvature $\kappa_{\Phi_A}$ with mobility $M_A$.

The result was discovered in an attempt to understand the meaning 
of a result of Sapiro and Tannenbaum [ST]. They
showed that for curves in the plane, motion by curvature to the
{\sl one-third} power, in both $C$ and $A(C)$,  is invariant under
area-preserving affine maps $A$.
But in the anisotropic case (especially when $\Phi$ is
non-differentiable), weighted mean 
curvature makes no reasonable sense under affine maps
unless the surface energy function is also transformed; 
the normal directions of $C$ can be quite
different from those of $A(C)$. But
once the surface energy is so transformed, then {\sl provided the
mobility is also transformed,} the motion which is affine invariant uses
$\kappa_\Phi$, not $\kappa_\Phi^{1/3}$.

As surveyed and extended in [OST], the [ST]
result is a special case of motions 
of the form where the time derivative in the normal direction
is equal to the second spatial derivative 
in a metric which is invariant under some group action.
The interpretation of this paper in those terms is that
one has a metric defined only on $C$, with distance along $C$ being
given by $\Phi$ and normal
to $C$ being given by $M$.  When one maps $C$, $\Phi$,
and $M$ by $A$, one
changes the metric in such a way that these ``distances'' remain invariant.

The original purpose of the Sapiro and Tannenbaum work, however, was
with regard to image processing: motion by curvature is useful
to remove noise.  In that case, there is no
physical Wulff shape, and the meaning of curvature as the rate of
decrease of length with area swept out is
irrelevant.  Conversely, it is not clear whether the result of this
paper has any useful implications for image processing.

A final note concerns the importance of the mobility.  That there
should {\sl be} a mobility factor is evident from considering
units: the units of velocity are distance over time, and the units
of curvature are one over distance. In our
case, even if mobility is initially 1 on each unit normal vector and so
``invisible," it {\sl must} be transformed and used in the $C_A$
motion, or the original motion is not preserved under affine maps. 
Mobility shows up naturally in other places as well.  For example, in
order for Wulff shapes to shrink homothetically under motion by
weighted mean curvature with mobility, $M$ must be proportional to
$\Phi$.  And in the anisotropic Allen-Cahn equation, the diffuse
interface approach to motion by weighted mean curvature, the sharp
interface limit naturally produces an extra factor in $M$ of $\Phi$
in addition to whatever diffuse mobility is initially put into the
Allen-Cahn equation.  

\section{Definitions}  For further details and references concerning
this section, see [T1] [T2].

2.1 The Wulff shape of $\Phi$ is the convex body
$W_\Phi = \{x: x\cdot p \le \Phi(p)\quad \forall p\};$
it is the equilibrium crystal shape.
If $\Phi$ is nondifferentiable in direction $n$, then
$W_\Phi$ has a facet with normal $n$; let $\Lambda_\Phi(n)$ be its length.

2.2.1 If $C$ is a twice differentiable curve at $x$ in $C$ and
$\Phi$ is twice differentiable at $n_C(x)$, then 
$\kappa_\Phi(C,x) = -T_C(x)\cdot {d\over ds} \nabla \Phi |_{n_C(x)}$
where ${d\over ds}$ denotes differentiation 
with respect to arc length along $C$.
That is, it
is minus the ``surface'' divergence of $\nabla \Phi$.  
One computes that if $C$ is
the graph of $x_2 = {1\over 2r} x_1 + o(x_1^2)$ near
(0,0), with oriented tangent (1,0) at (0,0), then
$\kappa_\Phi(C,(0,0))= - {1\over r}{\partial^2 \Phi \over
\partial p_1^2}|_{p=(0,-1)}$, the curvature times a weight.

2.2.2 If $C$ around $x$ is a line segment, of length $\ell$, and
$\Phi$ is nondifferentiable at $n_C(x)$,
then $\kappa_\Phi(C,x)
= -\sigma \Lambda_\Phi(n_C(x))
/\ell$, where $\sigma$ is 1 if $C$ is locally convex
at each end of the line segment, $\sigma$ is -1 if $C$ is locally
concave at each end, and $\sigma = 0$ if $C$ 
is convex near one end and is concave near the other.  

2.2.3 The definition of $\kappa_\Phi$ in the smooth case comes
from the fact that $\Phi$-first-variation linear operator on vector
fields $g$ is representable as $$-\int_{x\in C}  g(x)\cdot n_C(x) 
\kappa_\Phi(C,x)\, dx.$$  The non-differentiable case comes from
using non-local deformations.  If $n_C(x)$ does not exist, then either
there is an energy minimizing tangent cone to $C$ at $x$ or else
the weighted mean curvature there is infinite.

2.3 Note that the unit normal $n_C(x)$ is oriented here as the
exterior normal to a 
region in the plane locally bounded by $C$, as a special case of
the unit oriented normal $n_S(x)$ to a surface $S$ 
bounding a region in $R^d$. The curvature vector field $\vec\kappa$
points in the direction of maximum
{\sl decrease} of length under variations, and the sign of the
curvature here is such that $\vec\kappa = \kappa n_C$. 
The curvature of the boundary of a convex body is thus negative.

\section{Statements and Proofs}

Lemma 3.1. Given $v$ in $R^2$, if $q$ is normal to $L(v)$,
then $L^T q$ is normal to $v$.

Proof. $q\cdot L(v) = L^T q \cdot v$.
\medskip

Lemma 3.2. $W_{\Phi_A}= L(W_\Phi)$.

Proof. $W_{\Phi_A} = \{x: x\cdot q \le \Phi_A(q)\quad \forall q\}
= \{x: x\cdot q \le \Phi(L^T q)\quad\forall q\}
= \{x: x\cdot (L^T)^{-1}p \le \Phi(p)\quad \forall p\}
= \{x: L^{-1}x\cdot p \le \Phi(p)\quad \forall p\}
= \{Ly: y\cdot p \le \Phi(p)\quad\forall p\}
= L(W_\Phi).$
\medskip

Lemma 3.3. $\kappa_{\Phi_A}(A(C),A(x)) = \kappa_\Phi(C,x)$.

Proof.  
Consider first the case that $\Phi$ is not differentiable at
$n=n_C(x)$, so that the boundary of $W_\Phi$ has a 
facet (line segment) $W$ with normal $n$.  
Suppose $S$ is the maximal line segment in $C$ through $x$.  By
lemma 3.2, the length of $L(W)$ is $\Lambda_{\Phi_A}(q)$
for some $q$, and by lemma 3.1, $q =(L^T)^{-1}n/|(L^T)^{-1}n|$, which is
normal to $A(S)$.
Since $S$ is parallel to $W$, the ratio of
lengths of $L(W)$ and $A(S)$ is that of $W$ and $S$.
Thus $\kappa_\Phi(C,x) =\kappa_{\Phi_A}(A(C), A(x))$.

If $\Phi$ is differentiable at $n$, the lemma can be proved by
approximation by nondifferentiable $\Phi$. It can also be obtained by
the characterization of the weighted mean curvature as the rate of
decrease in surface energy with the volume (with $d$=2, area) swept
out, with a proof as in section 4 of this paper. Or it can be obtained
as an exercise in direct computation, as follows. Assume without loss
of generality that $n= (0,-1)$ and write $C$ locally as a graph $x_2 =
{1\over 2r}x_1^2 + o(x_1^2)$. One computes
$$-\kappa_{\Phi_A}(A(C),A(0,0)) = \Big(  div_{A(C)}
\nabla \Phi A|_{p=(x_1 ,x_1 ^2/(2r))}  \Big) \Big|_{x_1 =0}$$  
$$ ={L(0,-1) \over|L(0,-1)|}
\cdot  \Big|{d A(x_1 ,x_1 ^2/(2r)) \over dx_1 }\Big|^{-1}\Big({d\over dx_1 }
\big({d \Phi(L^T q) \over dq_1},{d \Phi(L^T q) \over
dq_2}\big)\Big|_ {q=n_{A(C)}(A(x_1 ,x_1 ^2/(2r)))} \Big)\Big|_{x_1 =0} $$
$$= {1\over r}{\partial^2 \Phi\over \partial p_1^2}\Big|_{p=(0,-1)},$$ 
remembering that  ${\partial^2 \Phi \over
\partial p_1 \partial p_2}(0,-1) = 0$ as a result of
$p\cdot \nabla \Phi(p) = \Phi(p)$ for all $p$, which in turn
follows from $\Phi(p) =|p| \Phi(p/|p|)$ for all nonzero $p$.

\medskip
Theorem 3.4. If $C(t)$ is a family of curves for $t \in [0,t_0]$,
with normal velocity 
$M\kappa_\Phi$ and if $C_A(t)$ is a family of curves with normal velocity
$M_A \kappa_{\Phi_A}$ and  $C_A(0) = A(C(0))$, then
$A(C(t) ) = C_A(t)$ for all $t\in [0,t_0]$. 

Proof.  By hypothesis, $A(C(t)) = C_A(t)$ when $t=0$; suppose it
remains true up until time $t \in [0,t_0)$.  
If $q$ is the unit
normal to $C_A(t) = A(C(t))$ at $A(x)$ for some $x$, then  $L^T q/|L^T
q|$ is the unit normal to $C(t)$ at $x$. If   
$ds$ is the incremental motion of $x$ in the direction $n_{C(t)}(x)$,
then $q\cdot L(L^T q/|L^T q|)\, ds  = |L^T q| ds$ is 
$A(ds)$, the normal component of the motion of $A(C(t))$ at $A(x)$. 
We now plug in $ds =M(L^T q/|L^T
q|)(-\sigma {\Lambda \over\ell}) dt$ 
to obtain  $A(ds) = M(L^T q)(-\sigma {\Lambda
\over\ell}) dt$.  But this is precisely the incremental normal motion of
$C_A(t)$, given Lemma 3 and the definition $M_A(q) = M(L^T q)$.  

Singularities can be present initially or can develop if $\Phi$ is
not even. In fact, one might be
dealing with two-dimensional grain growth where there are many
triple junctions at all times.  
Whatever way one handles such singular points for $\Phi$ (see [T2] for
example), they should then by definition be handled in the same way for
$\Phi_A$.  Since $A$ maps any such singular points in $C(t)$ to
similar singular points in $C_A(t)$, the proof above carries over.

\section{Higher dimensions} 

Suppose $S$ is a ($d$-1)-dimensional oriented
rectifiable surface in
$R^{d}$, $\Phi$ and $M$ are norms
on $R^{d}$, and $A$ is an affine map on $R^{d}$, with $A(x) =
L(x)+b$. Again define $\Phi_A(q)= \Phi(L^T q)$ and $M_A(q) =
M(L^T q)$.  The critical computation is
$$\Phi_A(A(S)) = \int_{y\in
A(S)} \Phi_A(n_{A(S)}(y)) dy $$
$$= \int_{y\in A(S)}{ \Phi(n_S(A^{-1}y)) \over
|(L^{-1})^T n_S(A^{-1}y)| }dy $$
$$ = \int_{x\in S} {\Phi(n_{S}(x))\over |(L^{-1})^T n_S(x)| }
J_2 L(n_S(x)) dx
$$
$$ = \int_{x\in S} \Phi(n_S(x)) |{\rm det} L| dx.$$
Since weighted mean curvature is the rate of decrease of surface energy
with volume, and the volume changes by a factor of $|{\rm det} L|$
under mapping by $A$, once again the weighted mean curvature of $S$
at $x$, with respect to surface energy $\Phi$, is the same as the
weighted mean curvature of $A(S)$ at $A(x)$, with respect to surface
energy $\Phi_A$.  The same proof as before shows that once the
mobility is also transformed, motion by weighted mean curvature with
mobility is affine invariant (i.e., $A(C(t)) = C_A(t)$).

If $\Phi$ is differentiable, then one can extend the results
further to the case of non-constant-coefficient
$\Phi$ and to diffeomorphisms $f: R^{d} \to R^{d}$ rather than
just affine maps $A$, since the preceding computation of surface energy
is entirely local. Given $\Phi$ depending on $x$ and $p$, one defines
$\Phi_f$ by $\Phi_f(y,q) = \Phi(x, (Df_x)^T q)$ where
$x=f^{-1}(y)$. Here
 $$\int_{y\in f(S)} \Phi_f(n_{f(S)}(y)) dy 
= \int_{x\in S} \Phi(x,n_S(x)) |{\rm det} Df_x| dx.$$
Again, weighted mean curvature stays the same, the mobility changes,
and motion by weighted mean curvature with mobility is affine
invariant. 

This final result has the following interpretation: not only can you
look at motion by weighted mean curvature with mobility at an angle,
but you can look at it through a
distorted lens. As long as you use your locally-measured $\Phi_f$ and
$M_f$, you will observe motion by weighted mean curvature with
mobility, provided the original motion is by the actual weighted mean
curvature with mobility.  

\section*{Acknowledgements}
This work was partially supported by a grant from the National Science
Foundation. The author is also grateful for the hospitality of MSRI,
where this paper was completed, and to Karen Almgren for a careful reading
of the paper.

\section*{References}

[OST] Peter J. Olver, Guillermo Sapiro, and Allen Tannenbaum,
Invariant geometric evolution of surfaces and volumetric smoothing, SIAM
J. Appl. Math. {\bf 57} (1997), 176-194.
\smallskip 

[ST]  Guillermo Sapiro and Allen Tannenbaum, On affine plane curvature
evolution, J. Functional Anal. {\bf 119} (1994), 79-120.
\smallskip 

[T1 Jean E. Taylor, Mean curvature and weighted mean curvature, Acta
metallurgica et materialia {\bf 40} (1992), 1475-1485.
\smallskip 

[T2] Jean E. Taylor, Motion of curves by crystalline curvature,
including triple junctions and boundary points, Differential Geometry,
Proceedings of Symposia in Pure Math. {\bf 51} (part 1) (1993), 417-438. 
\smallskip

[TCH]  Jean E. Taylor, John W. Cahn and Carol A. Handwerker, Geometric Models
of Crystal Growth,  Acta Metall. Mater. {\bf 40} (1992), 1443-1474.
\end{document}